\documentclass[a4paper,12pt]{article}

\usepackage{amsmath}
\usepackage{amssymb}
\usepackage{graphicx}
\usepackage{hyperref}
\usepackage{algorithm, algorithmicx, algpseudocode}

\newcommand{\bomega}{\boldsymbol{\omega}}
\newcommand{\J}{\mathrm{i}}
\newcommand{\E}{\mathrm{e}}

\graphicspath{{./Figures/}}

\title{Spherical polar coordinate transformation for integration of
  singular functions on tetrahedra}

\author{Michael Carley}

\begin{document}

%% ==================================%%
%% sample for unstructured abstract %%
%%==================================%%

\abstract{A method is presented for the evaluation of integrals on
  tetrahedra where the integrand has an integrable singularity at one
  vertex. The approach uses a transformation to spherical polar
  coordinates which explicitly eliminates the singularity and
  facilitates the evaluation of integration limits. The method is also
  implemented in an adaptive form which gives convergence to a
  required tolerance. Results from the method are compared to the
  output from an exact analytical method for one tetrahedron and show
  high accuracy. In particular, when the adaptive algorithm is used,
  highly accurate results are found for poorly conditioned tetrahedra
  which normally present difficulties for numerical quadrature
  techniques. The approach is also demonstrated for evaluation of the
  Biot-Savart integral on an unstructured mesh in combination with a
  fixed node quadrature rule and demonstrates good convergence and
  accuracy.}

\maketitle

%%\pacs[JEL Classification]{D8, H51}

%%\pacs[MSC Classification]{35A01, 65L10, 65L12, 65L20, 65L70}

\section{Introduction}
\label{sec:intro}

The integration of functions with an integrable singularity over
volume elements is a numerical operation which arises in a number of
fields, in particular where a volume potential is to be computed, such
as application of the Biot--Savart law in fluid dynamics and
electromagnetism, or in crack mechanics. In these cases, the
singularity in the integrand arises from the inverse distance
appearing in the Green's function for the problem, which has a
$1/R^{\alpha}$ dependence.

Owing to their importance in applications, a number of methods have
been developed over the years to evaluate such integrals, including
analytical~\cite{urankar84iv,weggel-schwartz88,onuki-wakao95,graglia87,%
  wilton-rao-glisson-schaubert-al-bundak-butler84,suh00,carley-angioni17},
semi-analytical~\cite{marshall-grant-gossler-huyer00}, and
fully-numerical~\cite{duffy82,khayat-wilton05,%
  lv-jiao-feng-wriggers-zhuang-rabczuk19} approaches. In this paper we
concentrate on the evaluation of volume integrals on tetrahedra, since
these often arise in applications using a structured or unstructured
mesh, and because they can be used to evaluate integrals on other
volume elements. The motivation for the work presented is the
evaluation of near-field terms in Fast Multipole Method (FMM)
accelerated application of the Biot--Savart law on volume meshes,
where far-field interactions can be handled using standard
quadratures, but a singular integration scheme is needed to correct
for near-field interactions, where the $1/R^{\alpha}$ singularity is
integrable, but is not well handled by standard quadrature rules.

Analytical approaches to the problem require an assumption about the
variation of source terms on the element. In these cases, the source
term is usually modelled as linear, as in methods which use the
divergence theorem to reduce the volume integral to a series of
line~\cite{suh00} or surface~\cite{marshall-grant-gossler-huyer00}
integrals, but monomial source terms of arbitrary order can also be
implemented~\cite{carley-angioni17}. For many purposes a fully
numerical approach is preferable, and a variable transformation is
often the most straightforward way to achieve this, as in the method
of Khayat and Wilton~\cite{khayat-wilton05} or the Duffy-type
transformations presented by Lv et
al.~\cite{lv-jiao-feng-wriggers-zhuang-rabczuk19}. 

The singularity considered in this paper is an inverse distance term,
which has been examined by various authors over many years. The
approach taken is typically to transform the variables of integration
so that the Jacobian of the transformation eliminates or smooths the
singularity in the integrand and allows standard one-dimensional
quadrature rules to be applied. The most straightforward of these is
similar to the methods used in dealing with surface integrals,
employing a transformation to cylindrical coordinates, accompanied by
a decomposition of the element into a number of sub-elements of a form
which facilitates the determination of integration limits. In the
Duffy transformation~\cite{duffy82}, analyzed and extended in a recent
paper~\cite{lv-jiao-feng-wriggers-zhuang-rabczuk19}, the tetrahedral
element is mapped onto a cube, making the Duffy transformation a
particular case of a ``pyramidal
transformation''~\cite{cano-moreno17}, an affine mapping of the
tetrahedron. In this case, the tetrahedron is decomposed into up to
three sub-elements of a form which allow the determination of
integration limits before application of the variable transformation.

In this paper, a method is presented which uses a transformation to
spherical polar coordinates. This approach appears to be novel and has
the advantage of explicitly eliminating the singularity without
requiring further variable transformations as in the Duffy
method~\cite{duffy82,lv-jiao-feng-wriggers-zhuang-rabczuk19}, and
allowing the use of standard one-dimensional Gaussian quadrature
rules. The only operation required is rotation of the tetrahedron to
an orientation which facilitates the evaluation of the integration
limits, Section~\ref{sec:integration:reference}, with the singularity
in the radial term being immediately removed by the change of
coordinate system. The second part of the method is the procedure for
rotating the tetrahedron into this reference orientation, in which the
integration limits can be easily calculated. This allows the method to
be applied to general tetrahedra, without requiring decomposition into
sub-elements, since the limits of integration are readily determined
for any tetrahedron in the reference orientation.

An estimate is presented of the convergence rate of the quadrature
method, which is confirmed by numerical testing using an analytical
formulation for integration on a tetrahedron. Further tests are
presented to demonstrate the performance of the adaptive method and of
the quadrature scheme used in a volume integration of the type which
appears in applications.

\section{Integration on tetrahedra}
\label{sec:integration}

The motivation for this paper is the evaluation of volume integrals on
tetrahedral meshes. It is assumed that most of the integration is
performed using standard fixed node quadrature rules on the elements.
Such rules are accurate for evaluation points far from an element, but
break down when the evaluation point lies on a tetrahedron. In this
case, we deal with integrals of the form
\begin{align}
  \label{equ:integral:basic}
  I &= \iiint_{V} \frac{f(\mathbf{x})}{R^{\alpha}}\,\mathrm{d}V,
\end{align}
on tetrahedra given by nodes $\mathbf{x}_{i}$, $i=0,1,2,3$, with the
singularity at node~0 so that
\begin{align}
  R &= \|\mathbf{x}-\mathbf{x}_{0}\|
\end{align}
where $\mathrm{d}V$ is the element of volume. In the applications
which motivate this work, the evaluation of volume potentials,
$\alpha=1,2$. Other values such as $\alpha=1/2$ arising in crack
mechanics, can be handled by a suitable choice of quadrature rule, and
results will also be presented for an irrational value of $\alpha$
close to~3, the limiting case for the integral to be integrable.

\subsection{Integration in the reference orientation}
\label{sec:integration:reference}

Integration is performed on the tetrahedron after transformation to a
reference orientation which facilitates the evaluation of integration
limits. In this orientation, the tetrahedron is defined by the
singular point, taken as the origin, and three nodes $\mathbf{y}_{i}$,
$i=1,2,3$. A spherical polar coordinate system centred at the origin
is used with
\begin{align}
  \label{equ:integration:sph}
  \mathbf{y} &= \rho
  \left[\sin\phi\cos\theta,\, \sin\phi\sin\theta,\, \cos\phi\right],
\end{align}  
noting that $\rho\equiv R$.

In the reference orientation, $\mathbf{y}_{1}=[\rho_{1},0,0]$ and
nodes $\mathbf{y}_{2}$ and $\mathbf{y}_{3}$ are given by
\begin{align}
  \label{equ:integration:nodes}
  \mathbf{y}_{2} &= \rho_{2}\left[\sin\phi_{2}\cos\theta_{23},\,
  \sin\phi_{2}\sin\theta_{23},\,\cos\phi_{2}\right],\\
  \mathbf{y}_{3} &= \rho_{3}\left[\sin\phi_{3}\cos\theta_{23},\,
  \sin\phi_{3}\sin\theta_{23},\,\cos\phi_{3}\right],
\end{align}
that is, the tetrahedron has been transformed so that node~1 lies on
the $(\theta,\phi)=(0,\pi/2)$ axis and nodes~2 and~3 lie in the
vertical plane given by $\theta=\theta_{23}$. The integral over the
tetrahedron is then
\begin{align}
  \label{equ:integration:sph:1}
  I &= \iiint \frac{f(\mathbf{x})}{R^{\alpha}}\,\mathrm{d}\mathbf{y}
  =
  \int_{0}^{\theta_{23}}
  \int_{\phi_{2}(\theta)}^{\phi_{3}(\theta)}
  \int_{0}^{\rho(\theta,\phi)}
  f(\mathbf{x})\rho^{n+\gamma}
  \mathrm{d}\rho\,
  \sin\phi\,
  \mathrm{d}\phi\,
  \mathrm{d}\theta,\\
  n + \gamma &= 2-\alpha,\nonumber
\end{align}
where $\phi_{2}(\theta)$ is the value of $\phi$ at which the vertical
plane with azimuthal angle $\theta$ intersects the line joining
nodes~1 and~2, and likewise for $\phi_{3}(\theta)$. In order to handle
non-integer values of $\alpha$, we introduce the notation
$n+\gamma=2-\alpha$, where $n$ is an integer and $-1\leq\gamma\leq1$.
Later, $\gamma$ will be used to define the weighting function of a
Gauss-Jacobi quadrature.

\begin{figure}
  \centering
  \includegraphics{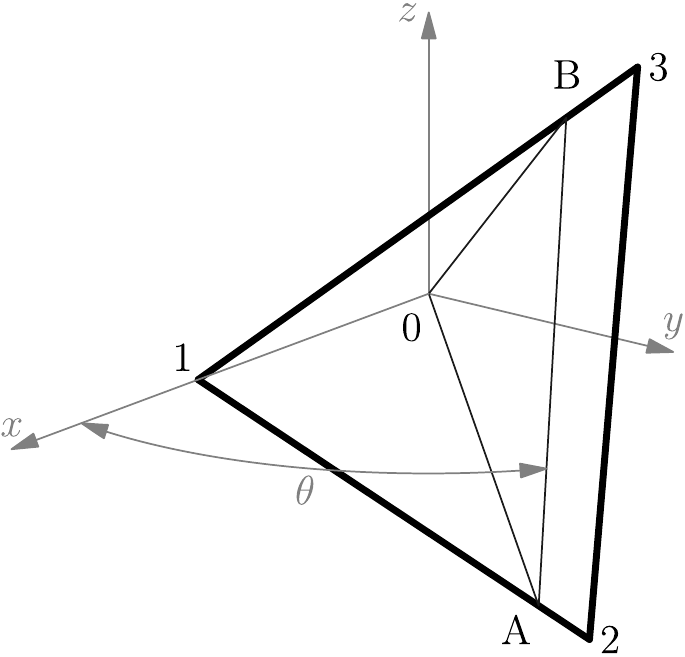}
  \caption{Integration on the tetrahedron in the reference orientation
    is performed over $\theta$, $0\leq\theta\leq\theta_{23}$. For each
    $\theta$, intersections of the $\theta_{23}$ plane with the
    face~123 (bold) are computed at the points A and B, and the limit
    on $\rho$ is computed on the line segment AB for each value of
    $\phi$}
  \label{fig:intersections}
\end{figure}

Integration over $\theta$ is performed between the fixed limits $0$
and $\theta_{23}$,
\begin{align*}
  \theta_{23} &= \tan^{-1}\frac{y_{2}}{x_{2}}
  = \tan^{-1}\frac{y_{3}}{x_{3}},
\end{align*}
with the node coordinates determined by the transformation procedure
of the next section.

The calculation of the limits in $\phi$ and $\rho$ requires two simple
calculations for the intersection of a plane with a line, and for the
intersection of two lines. Figure~\ref{fig:intersections} shows the
geometry of the system. At azimuthal angle $\theta$, the integration
limits in $\phi$ are determined by the intersection of the plane of
constant $\theta$ with the edges~12 and~13, denoted A and B
respectively. The intersection point is given by
\begin{align}
  \label{equ:intersect:2d}
  \mathbf{y}_{A} &= \mathbf{y}_{1} +
  (\mathbf{y}_{2}-\mathbf{y}_{1})u,\\
  % u &= -\frac{x_{1}\sin\theta - y_{1}\cos\theta}
  % {(x_{2}-x_{1})\sin\theta - (y_{2}-y_{1})\cos\theta},
  u &= -\frac{\mathbf{y}_{1}.\mathbf{s}}
  {(\mathbf{y}_{2}-\mathbf{y}_{1}).\mathbf{s}},\\
  \mathbf{s} &= [\sin\theta,\,-\cos\theta,\,0].
  \nonumber
\end{align}
In the spherical polar coordinate system, the limit $\phi_{A}$ is then
given by
\begin{align}
  \label{equ:intersect:2d:sph}
  \cos\phi_{A}
  &=
  \frac{u\rho_{2}\cos\phi_{2}}{\rho_{A}},\\
  \rho_{A}^{2} &=
  (1-u)^{2}\rho_{1}^{2}
  +u^{2}\rho_{2}^{2}
  +2u(1-u)\rho_{1}\rho_{2}\sin\phi_{2}\cos\theta_{23},\\
  u &= \frac{\rho_{1}\sin\theta}
  {\rho_{1}\sin\theta + \rho_{2}\sin\phi_{2}\sin(\theta_{23}-\theta)}.
\end{align}
A similar calculation is performed to find the limit $\phi_{B}$.
Integration over $\phi$ is then performed for
$\phi_{A}\leq\phi\leq\phi_{B}$ with $\phi_{A}$ and $\phi_{B}$ ordered
so that $\phi_{A}<\phi_{B}$.

For a given direction $(\theta,\phi)$, the limit on $\rho$ is
determined by the intersection of the ray through the origin with the
line segment AB. The distance from the origin to this intersection is
given by
\begin{align}
  \label{equ:intersect:rho}
  \rho &= \frac{\rho_{A}\rho_{B}\sin(\phi_{A}-\phi_{B})}%
  {(\rho_{B}\cos\phi_{B}-\rho_{A}\cos\phi_{A})\sin\phi - %
    (\rho_{B}\sin\phi_{B}-\rho_{A}\sin\phi_{A})\cos\phi}.
\end{align}
Integration on the tetrahedron can then be performed by integrating
over $\theta$, evaluating the limits of the inner integrals at each
point. 

\subsection{Transformation to reference orientation}
\label{sec:integration:transformation}

In the previous section a simple method of evaluating the singular
integral on a tetrahedron in a reference orientation was
presented. Given a point $\mathbf{y}$ on the tetrahedron, the
integrand $f(\mathbf{x})/\rho^{\alpha}$ can be evaluated using the
rotation matrix connecting the original and reference coordinate
systems and a translation of the origin,
\begin{align}
  \label{equ:integration:transform}
  \mathbf{x} &= \mathbf{x}_{0} + \mathbf{y}\mathsf{A},
\end{align}
where $\mathsf{A}$ is the rotation matrix.

The approach taken to transforming between coordinate systems is to
construct the transformed tetrahedron and then solve for
$\mathsf{A}$. Given the original tetrahedron nodes $\mathbf{x}_{i}$,
$i=0,1,2,3$, we define $\ell_{1}=\|\mathbf{x}_{1}-\mathbf{x}_{0}\|$
and set the first node of the transformed tetrahedron
$\mathbf{y}_{1}=[\ell_{1},0,0]$. Nodes~2 and~3 are positioned using
the constraint that the angle between line~01 and the plane containing
the triangle~023 must be the same in both coordinate systems. This is
achieved by the following procedure. First, calculate the normal to
the plane containing triangle~023, the projection $\mathbf{p}$ of
node~1 onto that plane, and construct a coordinate system centred at
$\mathbf{p}$. Writing $\mathbf{x}_{i}'=\mathbf{x}_{i}-\mathbf{x}_{0}$,
\begin{align}
  \label{equ:transform:1}
  \hat{\mathbf{n}} &= \frac{\mathbf{x}'_{3}\times\mathbf{x}'_{2}}
  {\|\mathbf{x}'_{3}\times\mathbf{x}'_{2}\|},\\
  d &= \mathbf{x}'_{1}.\hat{\mathbf{n}},\\
  \mathbf{p} &= \mathbf{x}'_{1} - d\hat{\mathbf{n}},\\
  \hat{\mathbf{s}} &= \mathbf{p}/\|\mathbf{p}\|,\\
  \hat{\mathbf{t}} &= \hat{\mathbf{n}}\times\hat{\mathbf{s}}.
\end{align}
This yields a coordinate system centred at point $\mathbf{p}$ with
unit vectors $\hat{\mathbf{s}}$, $\hat{\mathbf{t}}$, and
$\hat{\mathbf{n}}$, with the first two axes lying in the plane, and
the third being the normal to it. The angle between the edge~01 and
the plane~023 is then given by
\begin{align}
  \label{equ:transform:th}
  \theta_{23} = \cos^{-1}\frac{\mathbf{x}_{1}.\hat{\mathbf{s}}}{\ell_{1}}. 
\end{align}
To construct nodes~2 and~3 in the rotated coordinate system, we
establish a corresponding set of axes for the~023 plane as follows,
\begin{align}
  \label{equ:transform:10}
  \hat{\mathbf{n}}' &= \left[-\sin\theta_{23}, \cos\theta_{23}, 0\right],\\
  d' &= \mathbf{y}_{1}.\hat{\mathbf{n}}' = -\ell_{1}\sin\theta_{23},\\
  \mathbf{p}' &= \mathbf{y}_{1} - d'\hat{\mathbf{n}}',\\
  \hat{\mathbf{s}}' &= \mathbf{p}'/\|\mathbf{p}'\|,\\
  \hat{\mathbf{t}}' &= \hat{\mathbf{n}}'\times\hat{\mathbf{s}}'. 
\end{align}
This gives a coordinate system based on a plane with the correct
orientation with respect to $\mathbf{y}_{1}$ and allows the
calculation of nodes~2 and~3 as,
\begin{align}
  \mathbf{y}_{i} &=
  \mathbf{p}' +
  \hat{\mathbf{s}}.(\mathbf{x}'_{i}-\mathbf{p})\hat{\mathbf{s}}' + 
  \hat{\mathbf{t}}.(\mathbf{x}'_{i}-\mathbf{p})\hat{\mathbf{t}}'.
\end{align}
The rotation matrix $\mathsf{A}$ is then found by solving
\begin{align}
  \left[
    \begin{matrix}
      \mathbf{y}_{1}\\
      \mathbf{y}_{2}\\
      \mathbf{y}_{3}
    \end{matrix}
  \right]
  \mathsf{A}
  &=
  \left[
    \begin{matrix}
      \mathbf{x}'_{1}\\
      \mathbf{x}'_{2}\\
      \mathbf{x}'_{3}
    \end{matrix}
  \right],
\end{align}
for $\mathsf{A}$.

Algorithm~\ref{alg:integration} gives pseudocode for the evaluation of
the integral on a tetrahedron defined by four nodes $\mathbf{x}_{i}$,
$i=0,\ldots,3$, with the singularity at node~0. Required inputs are
Gauss-Legendre quadrature rules $(t_{i}^{(\theta)},w_{i}^{(\theta)})$,
$(t_{j}^{(\phi)},w_{j}^{(\phi)})$, for integration over $\theta$ and
$\phi$, respectively. Integration over $\rho$ is accomplished using a
Gauss-Jacobi rule $(t_{k}^{(\rho)},w_{k}^{(\rho)})$ with weight
function $(1+t)^{\gamma}(1-t)^{0}$. For integer values of $\alpha$,
$\gamma\equiv0$ and a Gauss-Legendre rule can be used. Rules are given
as $N_{\theta}$ nodes $t_{i}^{(\theta)}$ and weights
$w_{i}^{(\theta)}$, etc, with $-1<t_{i}^{(\theta)}<1$.

\begin{algorithm}
\caption{Pseudocode for integration on tetrahedron}\label{alg:integration}
\begin{algorithmic}
  \State generate tetrahedron in reference orientation
  $\mathbf{y}_{i}$ and solve for rotation matrix $\mathsf{A}$
  \State set $I=0$

  \State set $\bar{\theta}=\theta_{23}/2$,
  $\Delta\theta=\theta_{23}/2$

  \For{$i=1,\ldots,N_{\theta}$} calculate
  $\theta_{i}=\bar{\theta}+\Delta\theta t_{i}^{(\theta)}$

  \State calculate $\phi_{A}(\theta_{i})$, $\phi_{B}(\theta_{i})$ from
  Equation~\ref{equ:intersect:2d} or
  Equation~\ref{equ:intersect:2d:sph}. 

  \State set $\bar{\phi}=(\phi_{A}+\phi_{B})/2$,
  $\Delta\phi=\|\phi_{A}-\phi_{B}\|/2$

  \For{$j=1,\ldots,N_{\phi}$} calculate
  $\phi_{j}=\bar{\phi} + \Delta\phi t_{j}^{(\phi)}$

  \State calculate $\rho_{ij}(\theta_{i},\phi_{j})$ from
  Equation~\ref{equ:intersect:rho}

  \For{$k=1,\ldots,N_{\rho}$} calculate $\rho_{ijk}=\rho_{ij}(1 +
  t_{k}^{(\rho)})/2$

  \State set $\mathbf{y}_{ijk}=\rho_{ijk}
  [\cos\theta_{i}\sin\phi_{j},\sin\theta_{i}\sin\phi_{j},\cos\phi_{j}]$

  \State set $\mathbf{x}=\mathbf{x}_{0}+\mathbf{y}\mathsf{A}$

  \State set $\displaystyle I := I +
  f(\mathbf{x})\rho_{ijk}^{n}\sin\phi_{j}
  \Delta\theta\Delta\phi\left(\frac{\rho_{ij}}{2}\right)^{\gamma+1}
  w_{i}^{(\theta)}w_{j}^{(\phi)}w_{k}^{(\rho)}$
  \EndFor  

  \EndFor

  \EndFor
\end{algorithmic}
\end{algorithm}

\subsection{Adaptive quadrature}
\label{sec:integration:adaptive}

The basic algorithm of the previous section is easily implemented and
reliable for ``well-shaped'' tetrahedra. It does, however, have the
drawback that it gives no indication of the accuracy of the result
meaning that the user must either accept the possibility of
uncontrolled errors in the calculated value of the integral, or use
excessively high order quadrature rules with correspondingly excessive
computational effort. In this section, we give a method for using the
algorithm as the basis of an adaptive technique which can be used to
find results correct to some stated tolerance and which allows the
approach to be used with confidence on arbitrary tetrahedra.

\begin{figure}
  \centering
  \includegraphics{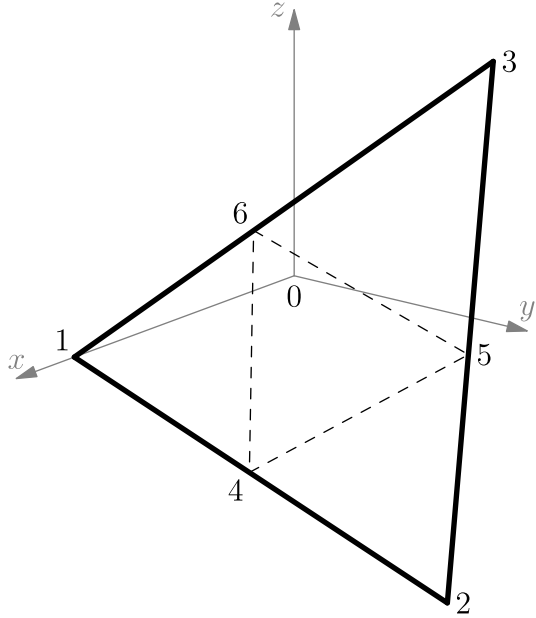}
  \caption{Splitting a tetrahedron for adaptive quadrature. The
    tetrahedron of Figure~\ref{fig:intersections} is shown with the
    base triangle split in four, generating four new tetrahedra each
    with a vertex at node~0.}
  \label{fig:integration:split}
\end{figure}

The procedure consists of splitting the face~123 into four new
triangles by bisecting each edge, as shown in
Figure~\ref{fig:integration:split}, and integrating over the resulting
sub-tetrahedra using the algorithm of the previous section. Given a
tolerance $\epsilon$, the convergence test is to check
\begin{align}
  \label{equ:integration:split}
  \left\|I_{0} - \sum_{i=1}^{4}I_{i}\right\| \leq \epsilon,
\end{align}
where $I_{0}$ is the integral evaluated on the original tetrahedron,
and $I_{i}$, $i=1,2,3,4$ is the integral on each
sub-tetrahedron. Sub-tetrahedra are~0146,~0425,~0536, and~0456. If the
convergence criterion is not met, the sub-tetrahedra are split and the
algorithm is applied recursively until the estimate of the total
converges to the required tolerance. The effectiveness of this
adaptive procedure will be demonstrated by numerical testing in
Section~\ref{sec:results:adaptive}.

\subsection{Convergence}
\label{sec:analysis:error}

The error behaviour of the integration method can be investigated by
considering an integral containing a polynomial source,
\begin{align}
  \label{equ:analysis:error:1}
  I &=
  \int_{0}^{\theta_{23}}
  \int_{\phi_{2}(\theta)}^{\phi_{3}(\theta)}
  \int_{0}^{\rho(\theta,\phi)}
  x^{i}y^{j}z^{k}\rho^{2-\alpha}
  \mathrm{d}\rho\,
  \sin\phi\,
  \mathrm{d}\phi\,
  \mathrm{d}\theta,\\
  N &= i+j+k,\nonumber
\end{align}
which in the spherical coordinate system is
\begin{align}
  \label{equ:analysis:error:2}
  I &=
  \int_{0}^{\theta_{23}}
  \int_{\phi_{2}(\theta)}^{\phi_{3}(\theta)}
  \int_{0}^{\rho(\theta,\phi)}
  \sin^{i+j+1}\phi\cos^{k}\phi\cos^{i}\theta\sin^{j}\theta
  \rho^{N+2-\alpha}
  \mathrm{d}\rho\,
  \mathrm{d}\phi\,
  \mathrm{d}\theta.
\end{align}
The error incurred by evaluating the integral numerically will depend
on the geometry of the tetrahedron and on the details of the three
quadrature rules used. To examine the error behaviour, we assume that
$\rho^{N+2-\alpha}$ can be integrated exactly. In this case, the error
depends on the trigonometric integrals in $\theta$ and $\phi$. We
begin by developing an error estimate for the evaluation of
\begin{align}
  \label{equ:analysis:error:3}
  I_{N} &= \int_{\phi_{1}}^{\phi_{2}}\sin^{N+1}\phi\,\mathrm{d}\phi,
\end{align}
using an $n$-point Gauss-Legendre quadrature. The quadrature rule
integrates exactly the expansion of the integrand in Legendre
polynomials up to order $2n-1$. Thus we derive an error estimate from
the expansion of the integrand in Legendre polynomials,
\begin{align}
  \label{equ:analysis:error:4}
  \sin^{N+1}\phi &= \sum_{i=0}^{\infty} a_{i}P(t),\quad -1\leq t \leq
  1,\\
  a_{i} &=
  \frac{2i+1}{2}\int_{-1}^{1}P_{i}(t)\sin^{N+1}\phi\,\mathrm{d}t,\nonumber\\
  \phi &= \bar{\phi} + t\Delta\phi,\nonumber\\
  \bar{\phi} &= \frac{\phi_{2}+\phi_{1}}{2},\quad
  \Delta\phi = \frac{\phi_{2}-\phi_{1}}{2},
  \nonumber
\end{align}
and $P_{i}(t)$ is the Legendre polynomial. The first term neglected in
the expansion of the integrand is $a_{2n}P_{2n}(t)$ and an upper bound
on the error is $2\|a_{2n}\|$. To estimate the coefficient $a_{2n}$,
we write
\begin{align}
  \label{equ:analysis:error:5}
  \sin^{N+1}\phi
  &=
  \frac{1}{(2\J)^{N+1}}
  \left(
    \E^{\J\phi} - \E^{-\J\phi}
  \right)^{N+1}\,\\
  &=
  \frac{\Delta\phi}{(\J 2)^{N+1}}
  \sum_{q}^{N+1}
  \binom{N+1}{q}
  (-1)^q
  \E^{\J(N+1-2q)\bar{\phi}}
  \E^{\J(N+1-2q)\Delta\phi t}.
\end{align}
Use of tables~\cite[7.324.2]{gradshteyn-ryzhik80} gives
\begin{align}
  \label{equ:analysis:error:bessel}
  \int_{-1}^{1}\E^{\J a t}P_{2n}(t)\,\mathrm{d}t
  &=
  (-1)^{n}\sqrt{\frac{2\pi}{\|a\|}}J_{2n+1/2}(\|a\|),
\end{align}
where $J_{\nu}$ is the Bessel function of the first kind. Using
Equation~\ref{equ:analysis:error:bessel}, the coefficient $a_{2n}$ is
\begin{align}
  a_{2n}
  =
  \Delta\phi\frac{4n+1}{2(2\J)^{N+1}}
  \sum_{q}^{N+1}
  &\binom{N+1}{q}
  (-1)^{q+n}
  \E^{\J(N+1-2q)\bar{\phi}}\times\nonumber\\
  &\sqrt{\frac{2\pi}{\|N+1-2q\|\Delta\phi}}J_{2n+1/2}
  \left(
    \|N+1-2q\|\Delta\phi
  \right).
  \label{equ:analysis:error:6}
\end{align}
To estimate the error bound, we use the large-order asymptotic form of
the Bessel function~\cite[10.19.1]{dlmf10},
\begin{align}
  J_{2n+1/2}(z)
  &\sim
  \frac{1}{\sqrt{2\pi(2n+1/2)}}
  \left(
    \frac{ez}{4n+1}
  \right)^{2n+1/2},
\end{align}
and hypothesize that for large $n$ the calculated integral should
converge as
\begin{align}
  \label{equ:analysis:error}
  \epsilon &\sim O(n^{-an})
\end{align}
for some constant $a$, which may depend on the geometry of the
tetrahedron.

\section{Results}
\label{sec:results}

Results are presented to demonstrate the performance of the quadrature
algorithm in computing integrals over tetrahedra of various
shapes. The first test integrals, which correspond to the evaluation
of a volume potential such as the Biot--Savart law, are
\begin{align}
  \label{equ:results:I}
  I_{ijk} &=
  \iiint
  \frac{x^{i}y^{j}z^{k}}{R^{\alpha}}\,
  \mathrm{d} x\,
  \mathrm{d} y\,
  \mathrm{d} z,\, i+j+k\leq N,
\end{align}
computed for $0\leq N\leq 4$, with $\alpha=1$. The relative error is
given as
\begin{align}
  \epsilon_{\mathrm{rel}} &=
  \max_{ijk}\frac{\left \| I_{ijk} - J_{ijk} \right \|}{\|J_{000}\|},
\end{align}
where $J_{ijk}$ is the exact value found using an analytical
method~\cite{carley-angioni17}. 

The geometries are chosen for comparison with previous work on a Duffy
transformation method~\cite{lv-jiao-feng-wriggers-zhuang-rabczuk19}
and test the algorithms on tetrahedra whose geometries cause
difficulties for numerical evaluation.

\subsection{Basic algorithm}
\label{sec:results:standard}

\begin{figure}
  \centering
  \includegraphics{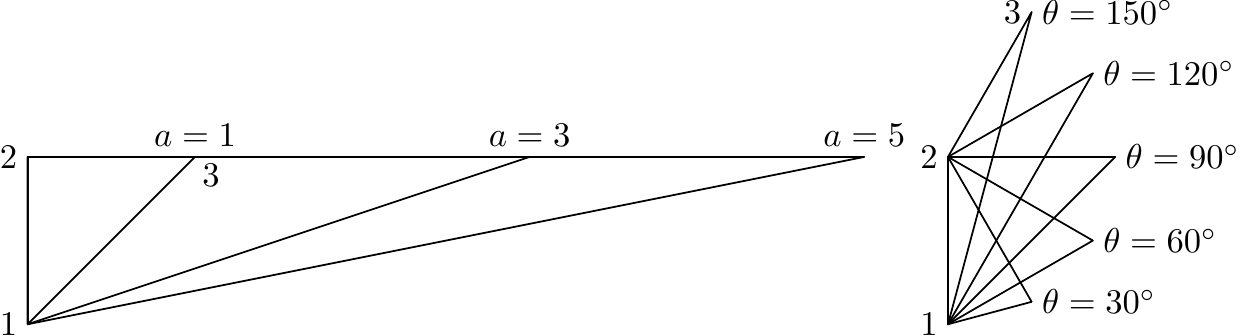}
  \caption{Variation of tetrahedron conditioning through changes in
    base triangle quality, with singular point $(0,0,0.1)$. Left hand
    side: increasing edge length $a$ increases the aspect ratio.
    Right hand side: changes in vertex angle with fixed edge
    lengths. Labels~1,~2,~3 correspond to the labelling of the base
    nodes in Figure~\ref{fig:intersections}.}
  \label{fig:results:geometry}
\end{figure}

We first present results to assess the performance of the basic
algorithm. The length of the Gauss-Legendre rules is varied, with
$N_{\theta}=N_{\phi}=N_{\rho}$ in each calculation. For compatibility
with previous work~\cite{lv-jiao-feng-wriggers-zhuang-rabczuk19},
three tetrahedra are considered. The first has nodes
$\mathbf{x}_{0}=[0,0,h]$, $\mathbf{x}_{1}=[0,0,0]$,
$\mathbf{x}_{2}=[0,1,0]$, and $\mathbf{x}_{3}=[1,1,0]$, with the
height $h$ varied to examine the effect on the quadrature error. The
second and third cases consider variations in the geometry of the base
of the tetrahedron, Figure~\ref{fig:results:geometry}, to study
possible effects of poorly-conditioned elements. In these cases, the
node $\mathbf{x}_{0}=[0,0,0.1]$. In the second case,
$\mathbf{x}_{3}=[\sin\theta, 1-\cos\theta, 0]$ to study the influence
of the vertex angle on the tetrahedron base. In the final test, the
effect of the base aspect ratio is considered, by setting
$\mathbf{x}_{3}=[a, 1, 0]$ and varying $a$.

\begin{figure}
  \centering
  \includegraphics{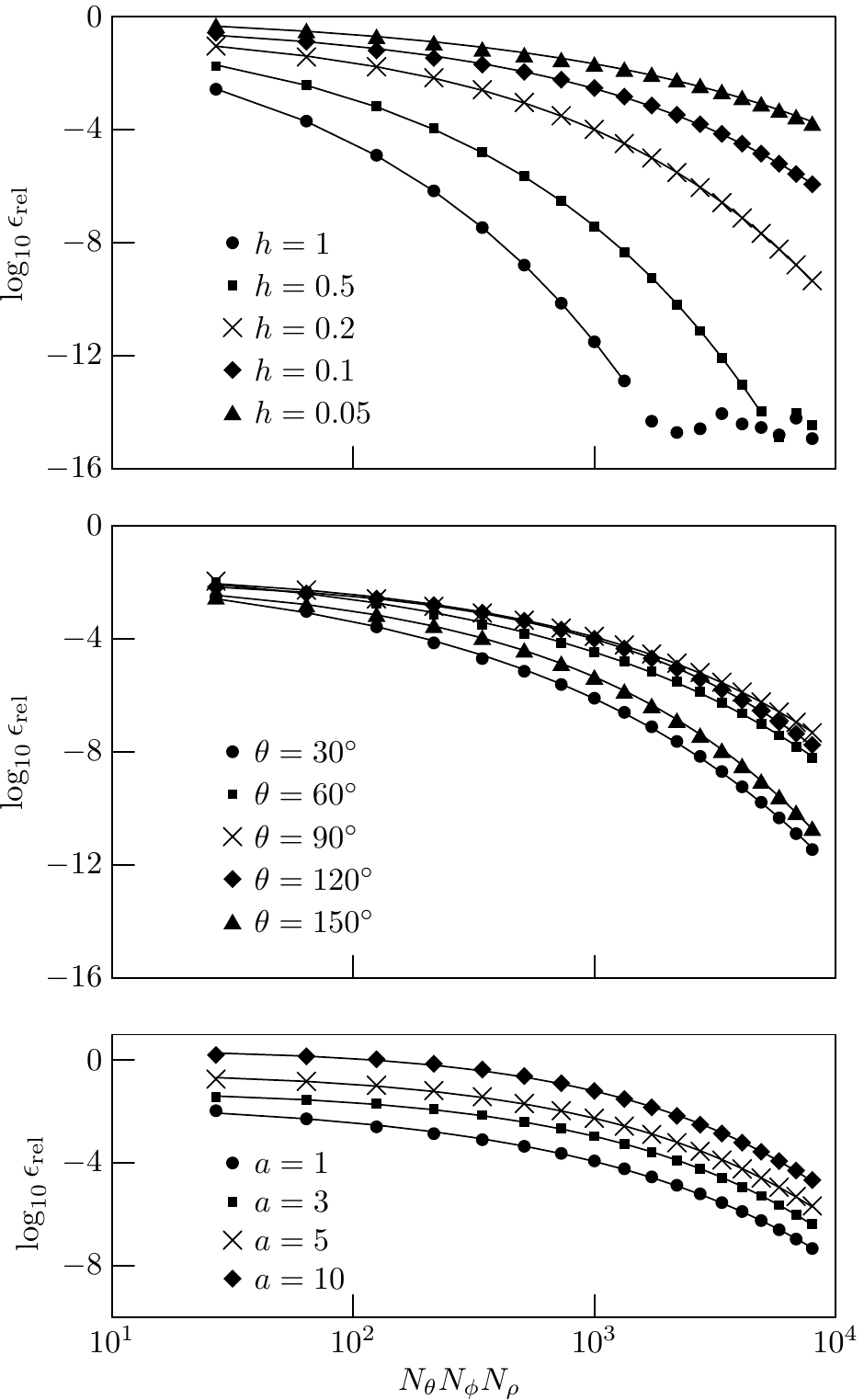}
  \caption{Effect of varying tetrahedron geometry on integration
    accuracy, error against total number of function evaluations; from
    top: varying tetrahedron height $h$, vertex angle $\theta$, base
    aspect ratio $a$. Symbols: error; solid lines: $An^{-an}$ fit.}
  \label{fig:results:basic}
\end{figure}

Figure~\ref{fig:results:basic} shows integration error on different
tetrahedra as a function of the number of quadrature points, including
the $n^{-an}$ convergence estimate, which is seen to fit the error
very well.

As might be expected, in the case of a ``well-conditioned''
tetrahedron, $h\approx1$ in the upper plot, convergence is rapid and
machine precision can be achieved. As $h$ is reduced, however, the
method cannot achieve high accuracy with the number of quadrature
points available and has quite poor results for $h=0.05$. In the other
two plots, the value of $h$ is held constant at~0.1 and the shape of
the tetrahedron base is varied. Again, as the tetrahedron becomes more
poorly conditioned, because of changes to the vertex angle or to the
base aspect ratio, the quadrature scheme shows poor accuracy and
unreliable convergence. The results of the next section will show how
the adaptive version of the algorithm can overcome these defects and
allow the method to converge to a required tolerance.

\begin{figure}
  \centering
  \includegraphics{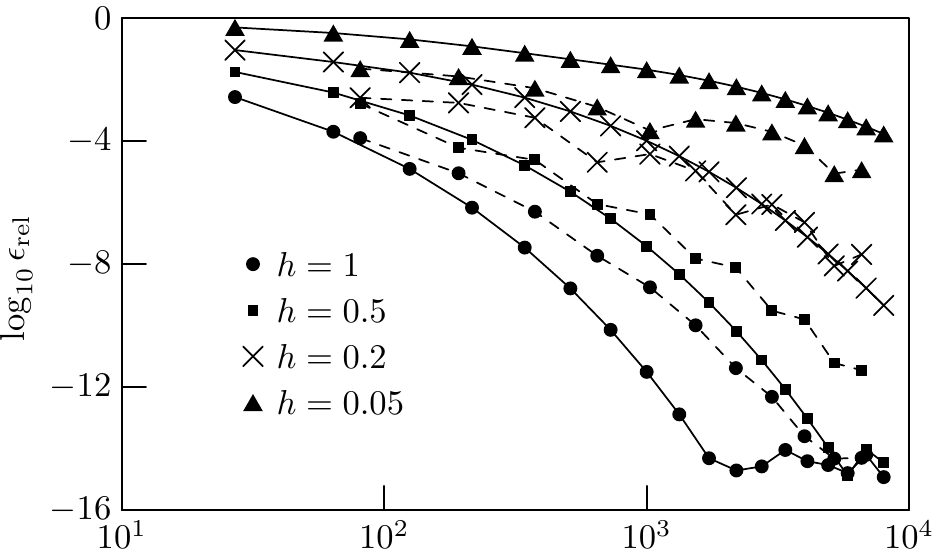}
  \caption{Error in test integral against number of function
    evaluations using the spherical coordinate method (solid lines)
    and the Duffy transformation (dashed), as a function of
    tetrahedron height.}
  \label{fig:results:duffy}
\end{figure}

With regard to computational effort, Figure~\ref{fig:results:duffy}
shows the number of function evaluations for integration on tetrahedra
of varying height, the same test case as in the first plot of
Figure~\ref{fig:results:basic}, using the spherical coordinate
transformation and a Duffy
transformation~\cite{lv-jiao-feng-wriggers-zhuang-rabczuk19}. For a
given quadrature rule length, the Duffy method has three times as many
evaluation nodes as the spherical coordinate method, because of the
decomposition of the tetrahedron in cylindrical
coordinates~\cite[page~15]{lv-jiao-feng-wriggers-zhuang-rabczuk19},
and this is accounted for in the operation count shown. For larger
$h$, i.e.\ better shaped tetrahedra, the spherical coordinate
transformation gives more rapid convergence, with the Duffy method
having better performance for $h\lessapprox0.1$, though the
convergence rate is quite poor in both cases for small $h$. 

Comparison of all of the test case data with the corresponding data
for a Duffy
method~\cite[Figures~18--22]{lv-jiao-feng-wriggers-zhuang-rabczuk19}
shows similar behaviour, though with slower convergence rates for the
spherical coordinate method as the shape of the base triangle is
changed. The convergence rate of the Duffy transform can be improved
using further changes of
variables~\cite{lv-jiao-feng-wriggers-zhuang-rabczuk19}; in this paper
we employ the adaptive approach of
Section~\ref{sec:integration:adaptive}, tested in the next
section. Overall, the computational effort for the Duffy-type methods
is comparable to that of the method presented here, with the total
number of integration points reaching values of the order of $10^{4}$
in order to achieve machine precision in some cases.

\subsection{Adaptive algorithm}
\label{sec:results:adaptive}

\begin{figure}
  \centering
  \includegraphics{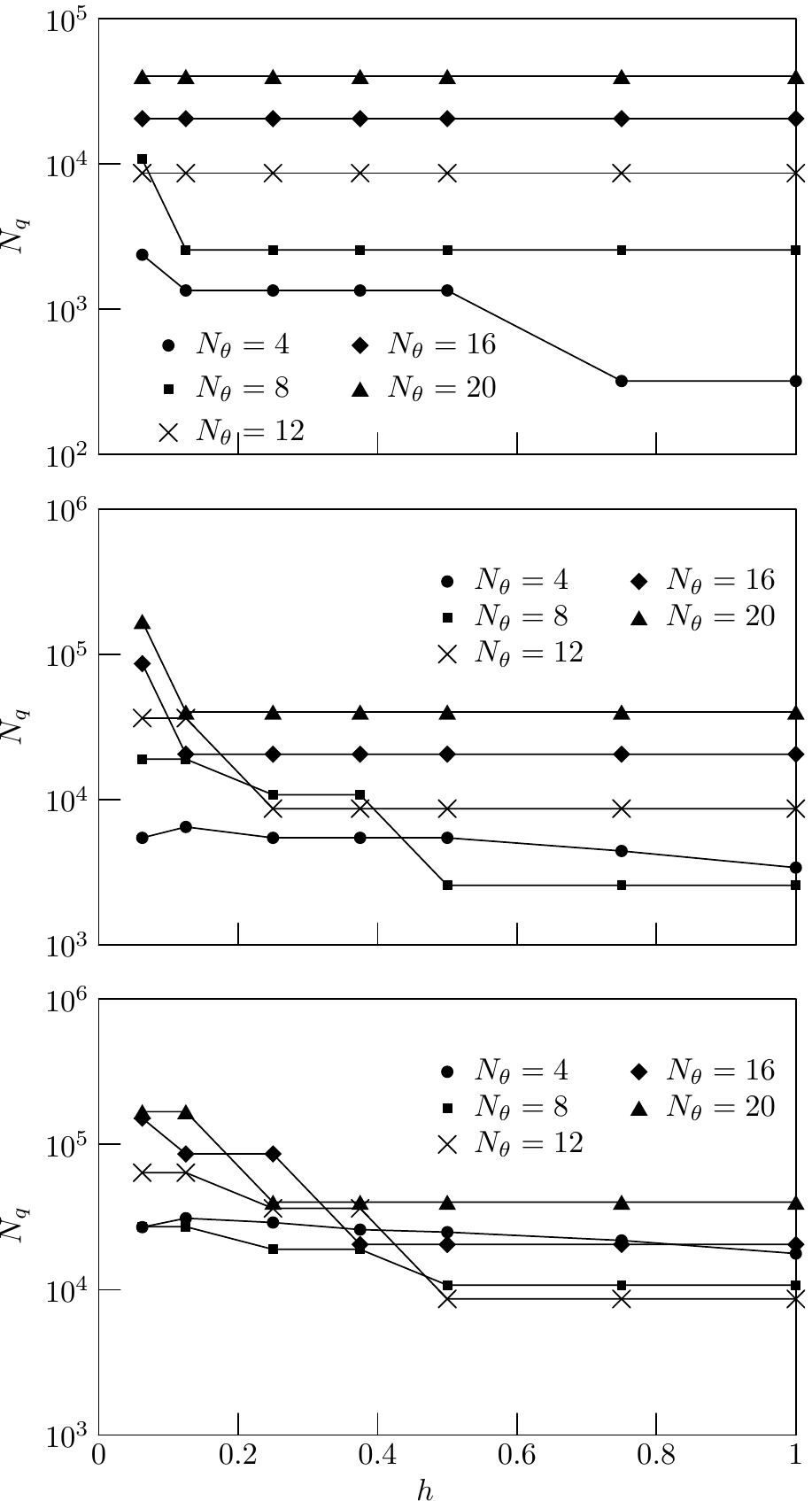}
  \caption{Effect of varying tetrahedron height on number of
    quadrature points at fixed tolerance using adaptive quadrature,
    total number of function evaluations against height $h$; from top:
    tolerance $\epsilon=10^{-3}, 10^{-6}, 10^{-9}$.}
  \label{fig:results:adaptive}
\end{figure}

Results for the adaptive version of the algorithm are presented in
Figure~\ref{fig:results:adaptive}. In this case the tetrahedron nodes
are $\mathbf{x}_{0}=[0,0,h]$, $\mathbf{x}_{1}=[0,0,0]$,
$\mathbf{x}_{2}=[0,1,0]$, and $\mathbf{x}_{3}=[2,1,0]$ with $h$ being
changed to modify the conditioning of the tetrahedron. As before, an
exact method is used to evaluate the integrals on the resulting
tetrahedra. The adaptive algorithm is applied for three tolerances,
$\epsilon=10^{-3},10^{-6},10^{-9}$ and using Gaussian quadrature rules
of length $N_{\theta}=4,8,\ldots,20$. Results presented are the total
number of function evaluations required to reach the requested
tolerance as a function of the tetrahedron conditioning represented by
the height $h$.

The results show the expected convergence behaviour. For $h\approx1$,
any of the quadrature rules gives a solution to the required
tolerance, as for the basic algorithm evaluated in the previous
section, but as $h$ is reduced and the tetrahedron becomes more poorly
conditioned, the method needs a greater number of recursions to
achieve convergence. At small values of $h$ and $\epsilon=10^{-9}$,
this leads to a large number of function evaluations. In all cases,
however, the requested tolerance is achieved, even when quite
low-order Gaussian quadratures, $N=4,8$, are used to integrate on
sub-elements. Numerical tests for accuracy when the base shape is
modified give similar results with convergence roughly independent of
aspect ratio $a$ and vertex angle $\theta$ and are not presented here.

It is interesting to note that the low-order rules can require more
function evaluations to achieve a given tolerance than higher order
rules, in particular, for $h\approx1$. This appears to happen because
for the well-conditioned tetrahedron, the low-order rules can require
more recursion levels to reach the convergence criterion.

We hypothesize that the adaptive algorithm achieves good convergence
because the base splitting shown in Figure~\ref{fig:integration:split}
has the effect of generating four tetrahedra which are better
conditioned than the parent element by virtue of having smaller area
bases with a constant height, in effect increasing $h$ with a
corresponding improvement in the element conditioning.

\subsection{Non-integer $\alpha$}
\label{sec:results:alpha}

The previous sections show how the proposed algorithm works for an
integrand with integer $\alpha$, where a reference value can be
evaluated exactly using analytical methods. This is an important case
in many applications, but there are also problems where non-integer
values appear. Here, we consider two cases, $\alpha=1/2$ which arises
in crack problems in solid mechanics, and an irrational value of
$\alpha$ which poses particular difficulties for coordinate
transformation schemes.

\begin{figure}
  \centering
  \includegraphics{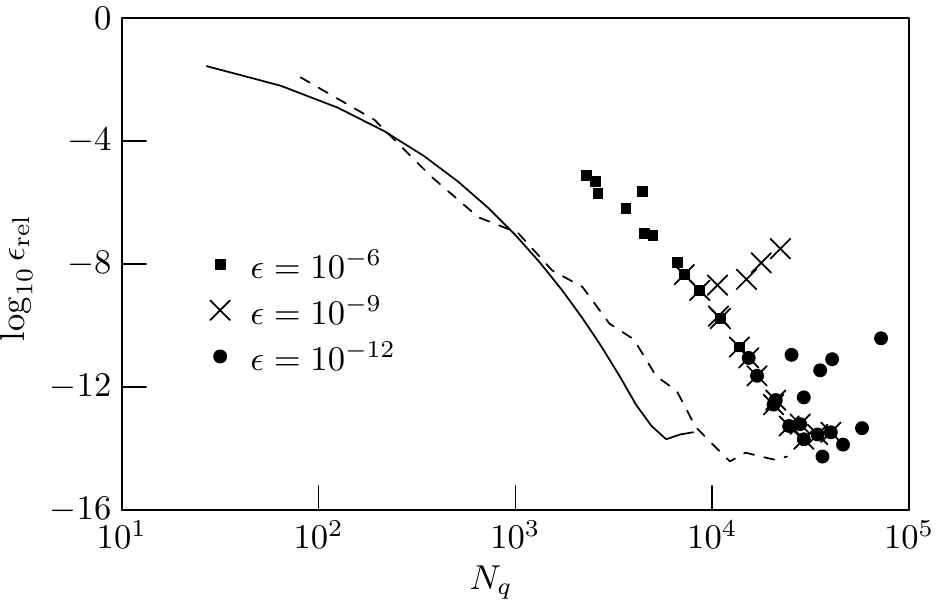}
  \caption{Convergence of algorithm for $\alpha=1/2$: error versus
    number of function evaluations. Solid line: basic algorithm;
    dashed line: Duffy transformation.}
  \label{fig:results:alpha:half}
\end{figure}

In the first case, that of rational $\alpha$, a Duffy-type
transformation~\cite{lv-jiao-feng-wriggers-zhuang-rabczuk19} can
eliminate the singularity and allow Gauss-Legendre rules to be used in
the radial direction. Figure~\ref{fig:results:alpha:half} shows the
performance of the method compared to the evaluation of a reference
integral using the Duffy-type
method~\cite{lv-jiao-feng-wriggers-zhuang-rabczuk19} with high order
Gaussian quadratures. The integral is evaluated on a tetrahedron with
nodes $\mathbf{x}_{0}=[0,0,1/2]$, $\mathbf{x}_{1}=[0,0,0]$,
$\mathbf{x}_{2}=[0,1,0]$, and $\mathbf{x}_{3}=[1,1,0]$. The Duffy
transformation was applied to the integral using high-order
Gauss-Legendre rules to give a reference value, and the method of this
paper was implemented using Gauss-Legendre rules in $\theta$ and
$\phi$ and a Gauss-Jacobi rule for $\rho$. In the non-adaptive case,
the computational effort for the spherical coordinate method is about
the same as for the Duffy transformation, though adaptive quadrature
incurs a computational cost to ensure convergence to the required
tolerance, as shown by the shift of the data points to the right.

\begin{figure}
  \centering
  \includegraphics{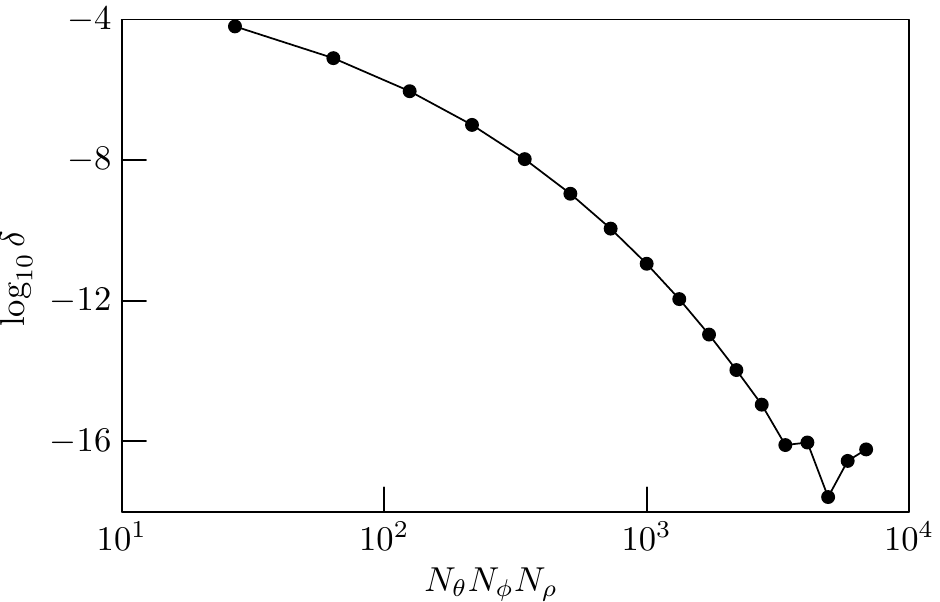}
  \caption{Convergence of algorithm for $\alpha=3-1/\pi$: difference
    between integral and reference versus number of function
    evaluations for basic algorithm.}
  \label{fig:results:alpha:3pi}
\end{figure}

Finally, we consider a case similar to that used in previous
work~\cite{chernov-von-petersdorff-schwab15}, with an irrational value
of $\alpha$ which is close to the point where the integrand is not
integrable. We set $\alpha=3-1/\pi$ and integrate using a quadrature
with $\gamma=1/\pi-1\approx-0.6817$. This integral cannot be evaluated
using the Duffy-type coordinate transformation, since it requires a
rational value of $\alpha$, so convergence is tested by evaluating
\begin{align}
  \delta &= \|I - I_{(\mathrm{ref})}\|,
\end{align}
where $I_{\mathrm{ref}}$ is the integral evaluated using quadrature
rules of length~20, or~$20^{3}=8000$ function evaluations, in the
spherical coordinate method. Figure~\ref{fig:results:alpha:3pi} shows
$\delta$ against the total number of function evaluations for the same
tetrahedron geometry as in Figure~\ref{fig:results:alpha:half}. In
order to track the convergence of the integral, a single monomial
source term is used, evaluating $I_{111}$ as defined in
Equation~\ref{equ:results:I}.

The convergence shown by Figure~\ref{fig:results:alpha:3pi} is quite
rapid, even for this demanding case, and is comparable to the
convergence shown in the corresponding plot in
Figure~\ref{fig:results:basic}. The Gauss-Jacobi rule handles the
singularity in the integrand and convergence to machine precision is
achieved.

\subsection{Evaluation of Biot--Savart integral}
\label{sec:results:biot}

The final numerical test evaluates the performance of the integration
method in a realistic application, the evaluation of the Biot--Savart
integral in three dimensions. This is an application which has
motivated the development of a number of integration
techniques~\cite[for
example]{urankar84iv,suh00,marshall-grant-gossler-huyer00} because of
its importance in fluid dynamics and electromagnetism. Here we
evaluate the velocity field of a vortex ring, a basic component of
many flows~\cite{shariff-leonard92}. The velocity induced by a
three-dimensional distribution of vorticity $\bomega(\mathbf{x})$ over
a volume $V$ is given by~\cite[page~18]{danaila-kaplanski-sazhin21}
\begin{align}
  \label{equ:results:biot}
  \mathbf{u}(\mathbf{x})
  &=
  -\frac{1}{4\pi}
  \int_{V}
  \frac{\mathbf{(\mathbf{x}-\mathbf{y})\times\bomega(\mathbf{y})}}{R^{3}}
  \,
  \mathrm{d} V,\\
  R &= \|\mathbf{x} - \mathbf{y}\|.\nonumber
\end{align}
For an axisymmetric ring, we employ cylindrical coordinates
$(r,\theta,z)$,
\begin{align}
  r^{2} &= x^{2} + y^{2},\\
  \theta &= \tan^{-1}\frac{y}{x}.
\end{align}
In the axisymmetric case, vorticity has only an azimuthal component
$\omega_{\theta}$ and the radial and axial velocity components are
given by~\cite[page~21]{danaila-kaplanski-sazhin21},
\begin{subequations}
  \label{equ:results:ring}
  \begin{align}
  u_{r}(r,z) &= -\frac{1}{r}\frac{\partial\psi}{\partial z},\\
  u_{z}(r,z) &= \frac{1}{r}\frac{\partial\psi}{\partial r},\\
  \psi &=
  \int_{r}\int_{z}
  \omega_{\theta}(r_{1},z_{1})
  \frac{\sqrt{rr_{1}}}{2\pi}
  \left[
    \left(
      \frac{2}{\kappa} - \kappa
    \right)K(\kappa)
    -\frac{2}{\kappa} E(\kappa)
  \right]
  \mathrm{d}r_{1}
  \mathrm{d}z_{1},\\
  \kappa^{2} &= \frac{4rr_{1}}{(z-z_{1})^{2} + (r+r_{1})^{2}},\nonumber
\end{align}
\end{subequations}
where $K(\kappa)$ and $E(\kappa)$ are complete elliptic integrals of
the first and second kind
respectively. Equations~\ref{equ:results:ring} can be used to evaluate
a reference velocity for comparison with the evaluation of
Equation~\ref{equ:results:biot}. As a test case, we use a
Gaussian-core ring with
\begin{align}
  \omega_{\theta} &=
  \exp
  \left[
    -\frac{(r-1)^{2} + z^{2}}{\sigma^{2}}
  \right],\\
  \bomega(\mathbf{x})
  &=
  \left(-\omega_{\theta}\sin\theta, \omega_{\theta}\cos\theta, 0\right),
\end{align}
and $\sigma=0.3$.

\begin{figure}
  \centering
  \includegraphics{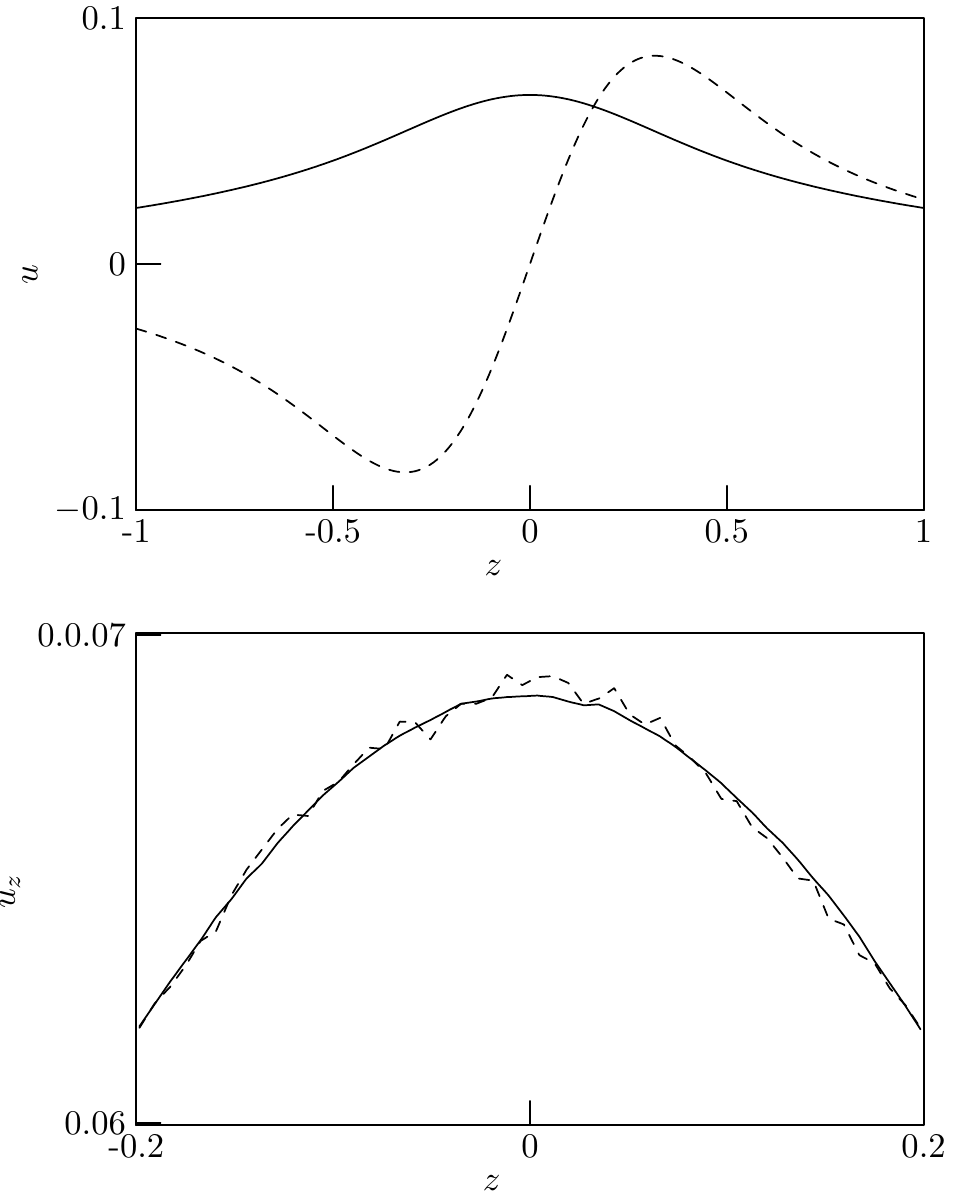}
  \caption{Axial (solid line) and radial (dashed line) velocity
    induced by Gaussian vortex ring at $r=1$ (top); zoom of axial
    velocity evaluated by integration on three-dimensional mesh with
    (solid line) and without (dashed line) singular integration
    method.}
  \label{fig:results:ring}
\end{figure}

For the three-dimensional evaluation, the vorticity is discretized on
an unstructured tetrahedral mesh with $-2.5\leq x\leq 2.5$,
$-2.5\leq y\leq 2.5$, $-1.5\leq z\leq 1.5$, using the \textsc{tetgen}
code~\cite{si15}, and the velocity field is evaluated at each node of
the tetrahedralization. Integration over tetrahedra is performed using
the high-order quadrature rules of Jaskowiec and
Sukumar~\cite{jaskowiec-sukumar21}, except for tetrahedra which have
an evaluation point as a vertex. In this case, the integral over the
tetrahedron is evaluated using the method of
Section~\ref{sec:integration}. For comparison with the axisymmetric
formulation, the velocity is evaluated at $(1, 0, z)$,
$-1\leq z\leq 1$, where
$\mathbf{u}(\mathbf{x})\equiv(u_{r}, 0,
u_{z})$. Figure~\ref{fig:results:ring} shows the radial and axial
velocity evaluated using the axisymmetric formulation. The lower plot
shows a zoom of the axial velocity evaluated using only fixed node
quadratures~\cite{jaskowiec-sukumar21} and the singular quadrature
method of this paper. The difference between the two plots is clear
and demonstrates the requirement for properly handling of the singular
integrand. 

Error in the calculation is controlled by the discretization of the
domain, by the order of the fixed node quadrature rules, and by the
order of the singular quadrature rules. The error measure used is
\begin{align}
  \label{equ:results:error:bs}
  \epsilon &=
  \frac{\max\|u'_{z}(z)-u_{z}(z)\|}{\max\|u_{z}(z)\|},\quad
  -1\leq z\leq 1,
\end{align}
where $u_{z}(z)$ is the axial velocity computed using
Equations~\ref{equ:results:ring}, and $u_{z}'(z)$ that computed using
the three-dimensional integration over tetrahedra. 

\begin{figure}
  \centering
  \includegraphics{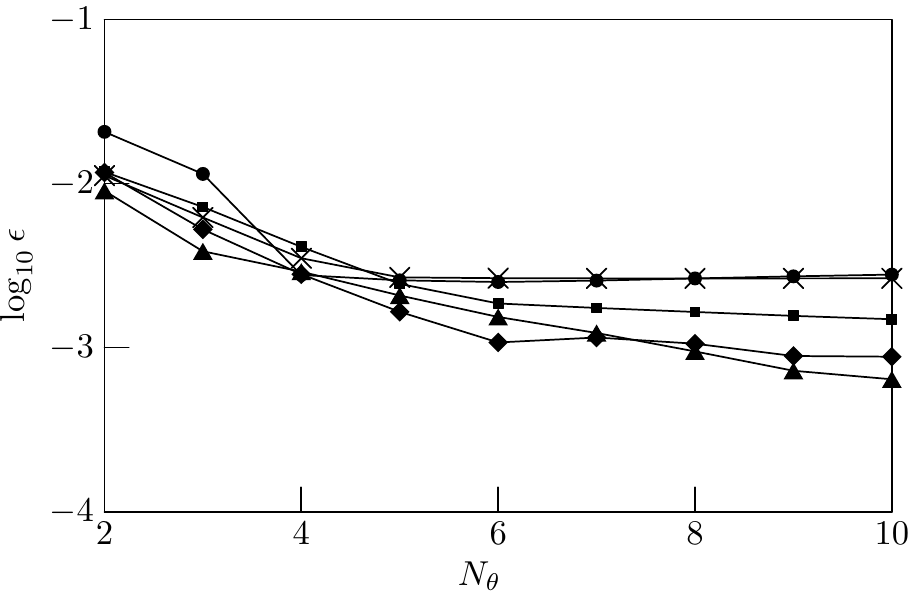}
  \caption{Error in axial velocity against number of singular
    quadrature points. Number of mesh points:~12626 (circle); 13260
    (box); 14117 (cross); 16712 (diamond); 35066 (triangle)} ;
  \label{fig:results:errors}
\end{figure}

Summation over the fixed quadrature points was performed using a fast
multipole method (FMM) based on the approach of Gumerov and
Duraiswami~\cite{gumerov-duraiswami03,gumerov-duraiswami04,
  gumerov-duraiswami05}. Direct summation was used at a set of sample
points to check that the truncation error of the FMM summation was at
least an order of magnitude less than the difference between the
computed and reference velocities.

The error in the evaluated velocity is given in
Figure~\ref{fig:results:errors}, which shows the error measure of
Equation~\ref{equ:results:error:bs} against the number of singular
integration points, as a function of the mesh discretization. For
these results, the twelfth order symmetric quadrature
rule~\cite{jaskowiec-sukumar21} was used for the non-singular
tetrahedron integration. As expected, meshes with a greater number of
nodes achieve a smaller error, though they require more quadrature
points to do so. The discretization error limit is reached for a
smaller number of quadrature nodes for coarser meshes than for the
more refined. 

\begin{table}[h]
\begin{center}
\begin{minipage}{174pt}
  \caption{Axial velocity error for~35066 node mesh}
  \label{tab:results:error}%
  \begin{tabular}{@{}rlllll@{}}
    \hline
    & \multicolumn{5}{c}{$N_{t}$}\\
    $N_{\theta}N_{\phi}N_{\rho}$& \multicolumn{1}{c}{11}
    & \multicolumn{1}{c}{23}
    & \multicolumn{1}{c}{44}
    & \multicolumn{1}{c}{74}
    & \multicolumn{1}{c}{117}
    \\
    \hline
    None & $2.5\times10^{-2}$& $3.4\times10^{-2}$& $4.1\times10^{-3}$&
    $1.9\times10^{-2}$& $2.8\times10^{-3}$\\ 
    8 & $6.5\times10^{-3}$& $5.3\times10^{-3}$& $5.3\times10^{-3}$&
    $5.1\times10^{-3}$& $5.1\times10^{-3}$\\ 
    27 & $7.4\times10^{-3}$& $3.3\times10^{-3}$& $1.8\times10^{-3}$&
    $1.9\times10^{-3}$& $1.8\times10^{-3}$\\ 
    64 & $8.2\times10^{-3}$& $3.1\times10^{-3}$& $9.8\times10^{-4}$&
    $1.4\times10^{-3}$& $8.5\times10^{-4}$\\ 
    125 & $8.2\times10^{-3}$& $3.0\times10^{-3}$& $1.1\times10^{-3}$&
    $1.4\times10^{-3}$& $8.5\times10^{-4}$\\ 
    216 & $8.0\times10^{-3}$& $2.9\times10^{-3}$& $1.1\times10^{-3}$&
    $1.1\times10^{-3}$& $6.1\times10^{-4}$\\ 
    343 & $7.8\times10^{-3}$& $2.9\times10^{-3}$& $1.1\times10^{-3}$&
    $8.4\times10^{-4}$& $5.0\times10^{-4}$\\ 
    512 & $7.7\times10^{-3}$& $2.9\times10^{-3}$& $1.1\times10^{-3}$&
    $7.3\times10^{-4}$& $4.9\times10^{-4}$\\ 
    729 & $7.7\times10^{-3}$& $2.9\times10^{-3}$& $1.1\times10^{-3}$&
    $7.5\times10^{-4}$& $4.8\times10^{-4}$\\ 
    \hline 
\end{tabular}
\end{minipage}
\end{center}
\end{table}

Table~\ref{tab:results:error} gives the error in evaluation of the
velocity on the finest mesh tested, as a function of the number of
singular and fixed quadrature points. In this case, the accuracy is
limited by the number of quadrature points The first row of the table
shows the error when the singular quadrature is not used and only the
non-singular method is applied. Lower rows show the effect of
including increasingly high-order singular rules and show that the
minimum error for a given non-singular quadrature is achieved when the
number of singular nodes $N_{\theta}N_{\phi}N_{\rho}$ is~2--3 times
the number of non-singular nodes $N_{t}$.

\section{Conclusions}
\label{sec:conclusions}

A method has been presented for the evaluation of integrals on a
tetrahedron with an integrable singularity at one vertex, motivated by
the evaluation of volume integrals used in fluid dynamics and
electromagnetism, and in fracture mechanics. The algorithm has been
shown to be reliable for well-conditioned tetrahedra in its basic
form. Extended to an adaptive form, it can compute the volume integral
to a required tolerance, even when the tetrahedron is poorly
conditioned. The method uses standard tools, such as one-dimensional
Gauss quadratures and simple geometric transformations, and can be
used without difficulty in production codes.

%% BioMed_Central_Bib_Style_v1.01

% \bibliography{abbrev,maths,vortex}% common bib file
%% if required, the content of .bbl file can be included here once bbl is generated
%%\input sn-article.bbl

%% Default %%
%%\input sn-sample-bib.tex%

\end{document}